# A Graph Invariant and 2-factorizations of a graph

Xie Yingtai

Chengdu University
(xtetai1@sina.com.cn)

**Abstract**

A spanning subgraph $\mathbb{R}$ of a graph G is called a [0,2]-factor of G, if $0 \leq d_{\mathbb{R}}(x) \leq 2$ for $\forall x \in V(G)$. $\mathbb{R}$ is a union of some disjoint cycles, paths and isolate vertices, that span the graph G. It is easy to get a [0,2]-factor of G and there would be many of [0,2]-factors for a G. A characteristic number for a [0,2]-factor, which reflect the number of the paths and isolate vertices in it,. The [0,2]-factor of G is called maximum if its characteristic number is minimum, and is called characteristic number of G.It to be proved that characteristic number of graph is a graph invariant and a polynomial time algorithm for computing a maximum [0,2]-factor of a graph G has been given in this paper.

A [0,2]-factor is Called a 2-factor , if its characteristic number is zero. That is ,a 2-factor is a set of some disjoint cycles, that span G. A polynomial time algorism for computing 2-factor from a [0,2]-factor, which can be got easily, is given..

A HAMILTON Cycle is a 2-factor, therefore a necessary condition of a HAMILTON Graph is that, the graph contains a 2-factor or the characteristic number of the graph is zero. The algorism, given in this paper, makes it possible to examine the condition in polynomial time.

**Key words:** 2-factor, [0,2]-factor ,Alternate Chain, P-chain , characteristic number of graph, graph invariant

## 1. Introduction

A 2-factor of Graph $G$ is a set of disjoint cycles that span $G$. 2-factors have multiple applications in Graph Theory, Computer Graphics, and Computational Geometry [1][2][3][4] A Hamiltonian cycle is then a 2-factor, and in one sense, it is the simplest 2-factor as it is composed of a single cycle. In another sense, it may be the most difficult 2-factor to find, as we must force a single cycle.

To our knowledge ,there are no efficient Algorithm for finding the 2-factor in general graph. In some cases an algorithm that computes such a 2-factor is also given. One such algorithm is due to Petersen. Petersen's result establishes the existence of 2-factors in $2m$-regular graphs only. Gopi and Epstein [5] propose an algorithm to compute 2-factors of 3-regular graphs. Their algorithm computes a perfect matching of the input graph. The edges, that is not in the computed matching define a 2-factor. Diaz-Gutierrez and Gopi [4] present two different methods to compute a 2-factors of graphs of maximum degree 4. The method consists of first computing a perfect matching on the input graph, after removing the edges in the matching from the input graph, and computing a new matching on the remaining subgraph. The 2-factor is defined by the edges in the union of both perfect matching. All appearance, their algorithm does not work on graphs with an odd number of vertices even when these are 4-regular. In fact, there are graphs with an even number of vertices where this algorithm also fails. The second method by Diaz-Gutierrez and Gopi is called the template substitution algorithm. In this method, vertices of degree 4 are replaced by templates, constructing by six vertex, to obtain an inflated graph. A perfect matching of the inflated graph can be translated into a 2-factor of the input graph given that exactly two outside vertices of a template connect to vertices of other templates. This algorithm is efficient for 4-regular only although the authors claim that their templates can also replace vertices of degree less than 4; however, no details are provided.

Umans [6] presented a known algorithm for computing 2-factors of a general graph. The algorithm is based on linear programming and, although simple to describe (as a set of linear equations), it is affected by the underlying complexity of linear programming ,but unfortunately, the problem about complexity of linear programming is yet a open problem[7][3].

We will propose a polynomial time algorithm for computing 2-factor of general graph. From any of a easily accessible [0,2]-factor to compute out the 2-factor in stages. This algorithm using a namely P-chain with respect to a [0,2]-factor of G is similar to augmenting path algogorism in the matching problem.

## 2. Basic theorem

The graph considered in this paper will be finite, undirected and simple. Let $G=(V,E)$ is a graph on $n$ vertices with vertex set $V(G)$ ($|V(G)|=n$) and edge set $E(G)$. Let S is a subgraph of G, $x \in V(G)$, the $d_S(x)$ is the degree of vertex $x$ in S. $P(x_1, x_2, ... x_k)$ is a path with two end vertices $x_1, x_k$ and inner vertices $x_i (1 < i < k)$, when $i \neq j$, $x_i \neq x_j$. $d_P(x_1) = 1, d_P(x_k) = 1$ and $d_P(x_i) = 2 (1 < i < k)$..

$C(x_1, x_2, ..... x_k, x_1)$ is a cycle with vertices $x_1, x_2, ... x_k$, $d_C(x_i) = 2 (i = 1, 2, ..., k)$ .when $k = n$, $P(x_1, x_2, ... x_n)$ is called HAMILTON path (H path), $C(x_1, x_2, ... x_n, x_1)$ is called HAMILTON cycle (H cycle).

Let $e = (x, y)$ is a edge with end vertices $x, y$ .Two edges $e_1, e_2$ is called adjacent if they have one (only one) common vertex. A sequence of adjacent edges without repeat $L = e_1 e_2 ... e_k$ is called a ***chain*** ( a ***close chain*** if $e_1, e_k$ is adjacent to each other, a ***open chain*** in the otherwise). Let $e_i = (x_i, x_{i+1})(i = 1, 2, ... k - 1)$ then $L(x_1, x_2, ... x_i, x_{i+1}, ... x_k)$ is a chain and $x_{i+1}$ is the common vertex of edges $e_i, e_{i+1}$. A path is also a chain, but a chain would not be necessarily a path,for example, in Fig.1(a), the $L(1, 5, 4, 8, 9, 6, 5, 12)$ is a chain only, but is not a path.

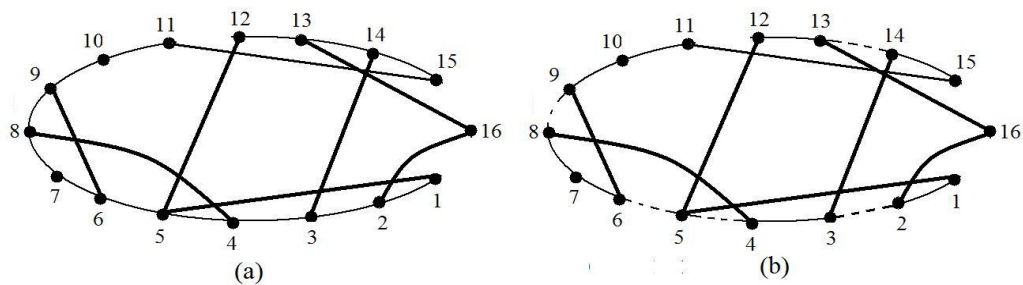

Fig 1.A [0,2]-factor and the P-chain with respect to it and the 2-factor of a Graph

**Definition 2.1** Let $\mathbb{R}$ is a spanning subgraph of graph G and $0 \leq d_\mathbb{R}(x) \leq 2$ for $\forall x \in V(G)$, we call $\mathbb{R}$ a [0,2]-factor of G.

$\mathbb{R}$ is a union of disjoint cycles , paths and isolated vertices that span G. Let

$$T_S(\mathbb{R}) = \sum_{v \in V(G)} (2 - d_\mathbb{R}(v)) \qquad (2.1)$$

Then $T_S(\mathbb{R})$ is a even and $T_S(\mathbb{R})/2$ is the number of paths and the isolated vertices in $\mathbb{R}$ .The $T_S(\mathbb{R})$ is called characteristic number of $\mathbb{R}$. When $T_S(\mathbb{R}) = 0$ $\mathbb{R}$ is 2-factor of G and If $\mathbb{R}$ is a single cycle then $\mathbb{R}$ is a H cycle. If G is a connected graph and $T_S(\mathbb{R}) = 2$ then there must be a H-path in G.

***Example 1*** In Fig.1,the two paths $P_1(1,2,3,4,5,6,7,8,9,10,11)$, $P_2(12,13,14,15)$ and a vertex $\{16\}$ in G form a [0,2]-factor, that $\mathbb{R} = P_1 \cup P_2 \cup \{16\}$ is a [0,2]-factor of G. And $T_S(\mathbb{R}) = 6$.

We here give a brief explanation for following work. First, It is easy to get a [0,2]-factor of G, for example ,the null graph of G(graph with V(G),but without any edges)is a namely [0,2]-factor. Let $\mathbb{R}_1$ is a [0,2]-factor of G, we will propose an Algorithm to get a [0,2]-factor $\mathbb{R}_2$ from $\mathbb{R}_1$, and $T_S(\mathbb{R}_2) = T_S(\mathbb{R}_1) - 2$,...,getting $\mathbb{R}_{i+1}$ from $\mathbb{R}_i$, and $T_S(\mathbb{R}_{i+1}) = T_S(\mathbb{R}_i) - 2$, until $\mathbb{R}_k$ with $T_S(\mathbb{R}_k) = 0$ is a 2-factor of G.

For this purpose, we need the following conception: we consider a graph G with a given partition of E(G) into disjoint subs $E_1$ and $E_2$, a chain $L = e_1 e_2 ... e_k$ is called *alternating,* if the adjacent edges in *L* are alternately in $E_1$ and $E_2$.

***Definition 2.2*** Let $\mathbb{R}$ is a [0,2]-factor of G. A chain $L = e_1 e_2 ... e_i e_{i+1} ... e_k$ is called P-chain with respect to $\mathbb{R}$, if *L* is a alternating respect to complementary subs E($\mathbb{R}$) and E(G) \ E(R),and satisfying following condition:

(1) $e_1, e_k \in E(G) \setminus E(\mathbb{R})$.

(2) when *L* is open then there are a end vertex *u* of $e_1$ and a end vertex *v* of $e_k$ such that $d_{\mathbb{R}}(u) \leq 1, d_{\mathbb{R}}(v) \leq 1$.

(3) when *L* is close then $e_1$ and $e_k$ with common end vertex *u* and $d_{\mathbb{R}}(u) = 0$.

Specially, we also consider the degenerate chain, when $k = 1$.

Intuitively speaking, a P-chain is a alternating chain, connecting two vertices $u, v$ which degree in $\mathbb{R}$ is less than two. and the two end edges do not belong to $E(\mathbb{R})$ ,that is ,to ask $e_1, e_k \in E(G) \setminus E(\mathbb{R})$.

The P-chain plays very important role in our algorithm to be established ,so we explain it more in details.

Let $e_i = (x_i, x_{i+1})(i = 1, 2, ... k)$ and $e_i, e_{i+1}$ with common vertex $x_{i+1}$ then $L(x_1, x_2, .... x_{k+1})$ is a P-chain with *k* edges $e_i$ ( $i = 1, 2, ... k$ ), if it is satisfying:

(1) *k* must be a odd. when *j* is odd $e_j \in E(G) \setminus E(\mathbb{R})$ ,when *j* is even $e_j \in E(R)$.

(2) if $x_1 \neq x_{k+1}$ then $d_{\mathbb{R}}(x_1) \leq 1, d_{\mathbb{R}}(x_{k+1}) \leq 1$ else if $x_1 = x_{k+1} = u$ then $d_{\mathbb{R}}(u) = 0$. $d_{\mathbb{R}}(x_j) = 2 (1 < j < k+1)$.

Because there are only two edges in $\mathbb{R}$ with common vertex $x_j$ for every $x_j(1 < j < k+1)$, so a P-chain *L* pass $x_j$ twice at most , therefor,

(3) $d_L(x_j) = 2, 4 (1 < j < k+1)$.

***Example 2 :*** In example 1(Fig.1), we have [0,2]-factor $\mathbb{R} = P_1 \cup P_2 \cup \{16\}$ ,then the alternate chains, $L_1(1,5,4,8,9,6,5,12), L_2(16,2,3,14,13,16)$ are the P-chain with respect to $\mathbb{R}$, $L_1$ is open and $L_2$ is close. Specially, the edge (11,15) is also a P-chain with respect to $\mathbb{R}$,that is degenerate.

We will to prove that if there is a [0,2]-factor $\mathbb{R}$ of G and a P-chain with respect to $\mathbb{R}$ then we can find a [0,2] -factor $\Re$ from $\mathbb{R}$ and the number of paths and isolated vertices in $\Re$ is less than the number of which in $\mathbb{R}$. This is what the Theorem 2.1 below claims.

**Theorem 2.1** Let $\mathbb{R}$ is a [0,2]-factor of G, $T_S(\mathbb{R}) > 0$, and $L(x_1, x_2, ... x_k, x_{k+1})$ ( $e_i = (x_i, x_{i+1})$ ) is a P-chain with respect to $\mathbb{R}$ then

$$\Re = \mathbb{R} \oplus L \ ^1$$

Is also a [0,2]-factor of G and

$$T_S(\Re) = T_S(\mathbb{R}) - 2$$

***Proof :*** First we prove that $0 \leq d_{\Re}(x) \leq 2$ for $\forall x \in V(G)$.

If $x \notin V(L)$ then $d_{\Re}(x) = d_{\mathbb{R}}(x)$ ,therefor, $0 \leq d_{\Re}(x) \leq 2$.

If $x \in V(L)$ then

---

[1] $G_1 \oplus G_2$ is a graph($G_1 \oplus G_2$)=V($G_1$) $\cup$ V($G_2$), E($G_1 \oplus G_2$)=E($G_1$) $\oplus$ E($G_2$) (A $\oplus$ B=(A $\cup$ B) \ (A $\cap$ B))

Case 1): We first consider the two end vertices $x_1, x_{k+1}$ of $L$, it can be easily seen from definition 2.2 that $d_{\Re}(x_1) = d_{\mathbb{R}}(x_1) + 1$ and $d_{\Re}(x_{k+1}) = d_{\mathbb{R}}(x_{k+1}) + 1$ when $x_1 \neq x_{k+1}$ (For example the vertices 1 and 12 in Fig .1) and $d_{\Re}(x_1) = d_{\Re}(x_{k+1}) = 2$ when $x_1 = x_{k+1}$ (For example the vertex 16 in Fig. 1).

Case 2): $x = x_i (i \neq 1, k+1)$ then $d_L(x_i) = 2, 4$

$2_a$) If $d_L(x_i) = 2$ and $d_{\mathbb{R}}(x_i) = 2$ then there are three edges with common vertex $x_i$. One and only one among the three edges belongs to $E(L) \cap E(\mathbb{R})$, threfor $d_{\Re}(x_i) = 2$ (for example the vertices 4,8,9,6 in Fig.1).

$2_b$) If $d_L(x_i) = 4$ and $d_{\mathbb{R}}(x_i) = 2$ then there are four edges with common vertex $x_i$. Two and only two among the four edge belongs to $E(L) \cap E(\mathbb{R})$, therefor $d_{\Re}(x_i) = 2$ (for example the vertices 5 in Fig.1).

we have proved that $d_{\mathbb{R}}(x) = d_{\Re}(x)$ for $\forall x \in V(G)$ except $x_1$ and $x_{k+1}$. And or $d_{\Re}(x_1) = d_{\Re}(x_{k+1}) = 2$ (when $x_1 = x_{k+1}$, in this case $d_{\mathbb{R}}(x_1) = d_{\mathbb{R}}(x_{k+1}) = 0$) or $d_{\Re}(x_1) = d_{\mathbb{R}}(x_1) + 1$ and $d_{\Re}(x_{k+1}) = d_{\mathbb{R}}(x_{k+1}) + 1$ (when $x_1 \neq x_{k+1}$), as a result of case 1 in proof. Hence one can to conclude that: $\Re$ is also a [0,2]-factor of G and $T_S(\Re) = T_S(\mathbb{R}) - 2$. Theorem has been proved.

□

***Example 3.*** In fig.1, by example 1, example 2, let $\mathbb{R}_1 = \mathbb{R} \oplus L_1$, then

$$\mathbb{R}_1 = P(15, 14, 13, 12, 5, 1, 2, 3, 4, 8, 7, 6, 9, 10, 11) \cup \{16\}$$

is also a [0,2]-factor of G. And $T_S(\mathbb{R}_1) = T_S(\mathbb{R}) - 2 = 4$. And the edge (11,15) is a P-chain with respect to $\mathbb{R}_1$, therefor the [0,2]-factor of G

$$\mathbb{R}_2 = \mathbb{R}_1 \oplus (11,15) = C(15, 14, 13, 12, 5, 1, 2, 3, 8, 7, 6, 9, 10, 11, 15) \cup \{16\}$$

Formed by a cycle and a isolated vertex and $T_S(\mathbb{R}_2) = 2$. And the chain $L_2(16, 2, 3, 14, 13, 16)$ is also a P-chain with respect to $\mathbb{R}_2$, then

$$\mathbb{R}_3 = \mathbb{R}_2 \oplus L_2 = C_1(1, 5, 12, 13, 16, 2, 1) \cup C_2(3, 14, 15, 11, 10, 9, 6, 7, 8, 4, 3)$$

Is a 2-factor of G, as shown in fig.1 (b).

**Definition 2.3** The [0,2]-factor $\mathbb{R}$ of G is ***maximum***, if which characteristic number $T_S(\mathbb{R})$ is ***minimum.*** The characteristic number of a maximum [0.2]-factor is called characteristic number of Graph G.

There can be no P-chain with respect a maximum [0,2]-factor since such a P-chain can be used, by Theorem 2.1, to get a [0,2]-factor with smaller Characteristic number. It turns out that the converse is true as well.

**Theorem 2.2** A [0,2]-factor $\mathbb{R}$ of G is maximum if and only if there is no P-chain with respect to $\mathbb{R}$.

*Proof* : One direction follows from Theorem 2.1. For other direction, we suppose that there is no P-chain in G with respect to $\mathbb{R}$, and yet $\mathbb{R}$ is not maximum. That is, there is a [0,2]-factor $\Re$ of G such that $T_S(\Re) < T_S(\mathbb{R})$, we consider the edges in $\mathbb{R} \oplus \Re$; these edges together with their end vertices form a subgraph $\mathbb{Q} = (V, \mathbb{R} \oplus \Re)$ of G (this subgraph may be disconnected). We will to prove that there must be a P-chain with respect to $\mathbb{R}$ in $\mathbb{Q}$ (thus, in G). The main idea of proving theorem 2.2 is coming from the proof of similar theorem in theory of maximum matchings, but here is more complicated and more difficult.

First, we have

$$T_S(\mathbb{R}) > T_S(\Re) \Rightarrow \sum_{v \in V(G)} d_{\mathbb{R}}(v) < \sum_{v \in V(G)} d_{\Re}(v) \Rightarrow |E(\mathbb{R})| < |E(\Re)|$$

Thus

$$|E(\mathbb{R}) - E(\mathbb{R}) \cap E(\Re)| < |E(\Re) - E(\mathbb{R}) \cap E(\Re)| \qquad (2.2)$$

Let
$$E_{\mathbb{Q}}(\mathbb{R}) = E(\mathbb{R}) \setminus E(\mathbb{R}) \cap E(\mathfrak{R})$$
$$E_{\mathbb{Q}}(\mathfrak{R}) = E(\mathfrak{R}) \setminus E(\mathbb{R}) \cap E(\mathfrak{R})$$

then E($\mathbb{Q}$) is parted into two disjoint subs $E_{\mathbb{Q}}(\mathbb{R})$ and $E_{\mathbb{Q}}(\mathfrak{R})$. Above inequation shows that $|E_{\mathbb{Q}}(\mathbb{R})| < |E_{\mathbb{Q}}(\mathfrak{R})|$. This is a key fact to be used for our proof.

Those vertices $u$ in $V(\mathbb{Q})$ are called **_terminus_** if $d_{\mathbb{R}}(u) \leq 1$ or $d_{\mathfrak{R}}(u) \leq 1$. Let $u$ is a terminus and $d_{\mathfrak{R}}(u) \leq 1$, then the edge $e = (u,v) \in E_{\mathbb{Q}}(\mathbb{R})$ with $u$ as a end vertex is called $\mathbb{R}$-**_end-edge_**. Similar, a $\mathfrak{R}$-**_end-edge_** $e$ is such a edge that $e \in E_{\mathbb{Q}}(\mathfrak{R})$ with $u$ as a end vertex and $d_{\mathbb{R}}(u) \leq 1$. A **_end-edge_** either is a $\mathbb{R}$-**_end-edge_** or is a $\mathfrak{R}$-**_end-edge_**.

Those vertices $u$ in $V(\mathbb{Q})$ are called **_inner_** vertices if its degree is 2 in both $\mathbb{R}$ and $\mathfrak{R}$ (i.e. $d_{\mathbb{R}}(u) = 2$ and $d_{\mathfrak{R}}(u) = 2$). Obviously, if $u$ is a inner vertex then or $d_{\mathbb{Q}}(u) = 2$ and there is a twain alternate edge with common vertex $u$ or $d_{\mathbb{Q}}(u) = 4$ and there are two twain alternate edges with common vertex $u$. Another kind of vertices $u$ are such that $d_{\mathbb{Q}}(u) = 3$, one can descry that in this case $u$ is a terminus for a end-edge and is also a inner vertex for a twain alternate edges.

Let $L = e_1 e_2 \ldots e_k$ is a alternate chain with respect to $E_{\mathbb{Q}}(\mathbb{R})$ and $E_{\mathbb{Q}}(\mathfrak{R})$ in $\mathbb{Q}$ with two end-edges $e_1$ and $e_k$ (We also consider the degenerate, that is $k=1$).

Now to prove that $\mathbb{Q}$ can be decomposed[†] into above alternate chains of $\mathbb{Q}$. In fact, let $x_1$ is a terminus (It must exist in $V(\mathbb{Q})$ because $T_S(\mathfrak{R}) < T_S(\mathbb{R})$ )then there is a end-edge $e_1 = (x_1, x_2)$, if $x_2$ is also a terminus then $e_1$ is already a alternate chain(degenerate), if $x_2$ is a inner vertex then there is a alternate (with $e_1$) edge $e_2 = (x_2, x_3)$, preceding in this fashion, we can get a alternate chain $L_1 = e_1 e_2 \ldots e_k$ with two end-edges $e_1, e_k$. we consider the edges in $\mathbb{Q} \oplus L_1$; these edges together with their end vertices form a subgraph $\mathbb{Q}_1 = (V_1, \mathbb{Q} \oplus L_1)$ of G (this subgraph may be disconnected). It can be easily validated that every terminus in $\mathbb{Q}_1$ is also a terminus in $\mathbb{Q}$ and every inner vertex in $\mathbb{Q}_1$ is also a inner vertex in $\mathbb{Q}$. Thus we can ditto get a alternate chain $L_2 = e_{2_1} e_{2_2} \ldots e_{2_k}$ with two end-edge $e_{2_1}, e_{2_k}$ in $\mathbb{Q}_1$ (but is also in $\mathbb{Q}$). Proceeding in this fashion, we can get: $\mathbb{Q} = \mathbb{Q}_0, \mathbb{Q}_1, \mathbb{Q}_2, \ldots \mathbb{Q}_m = \phi$ and $L_1, L_2, \ldots L_m$ such that
$$\mathbb{Q}_i = (V_i, \mathbb{Q}_{i-1} \oplus L_i)$$

And $L_{i+1}$ is a alternating chain in $\mathbb{Q}_i$ (Thus, is also in $\mathbb{Q}$). We have already decomposed $\mathbb{Q}$ into $m$ alternate chains with two end-edge, by this way. A alternate chain $L$ with two end-edges have only one among following three form:

(1) Two end-edge of $L$ are $\mathfrak{R}$-end-edge.
(2) Two end-edge of $L$ are $\mathbb{R}$-end-edge.
(3) One of two end-edges of $L$ is $\mathfrak{R}$-end-edge and another is $\mathbb{R}$-end-edge.

Let $L_\lambda$ to denote the length of $L$, $L_{\mathbb{R}}$ denote the number of edges belonging to $E_{\mathbb{Q}}(\mathbb{R})$ and $L_{\mathfrak{R}}$ denote the number of edges belonging to $E_{\mathbb{Q}}(\mathfrak{R})$. In form (1) the $L_\lambda$ must be a odd and $L_{\mathbb{R}} < L_{\mathfrak{R}}$, in form (2) the $L_\lambda$ is also a odd but $L_{\mathbb{R}} > L_{\mathfrak{R}}$, in form(3) the $L_\lambda$ is even and $L_{\mathbb{R}} = L_{\mathfrak{R}}$. One can come to a conclusion that there must be a alternate chain in $\mathbb{Q}$ formed as (1) because (2.2). A alternate chain formed as (1) is just a P-chain with respect to $\mathbb{R}$. However, this contradicts our assumption that there is no P-chain in G with respect to $\mathbb{R}$. The theorem has been proved. □

Follows from above two theorem, one can draw conclusion that:

**Basic Theorem** The Characteristic number of a Graph is a Graph invariant, and a 2-factor in a graph G exists if and only if the Characteristic number of the graph is zero.

---

[†] The graph G is decomposed into subgraph $G_1, G_2$ if $G_1 \cup G_2 = G$ and $G_1 \cap G_2 =$ null graph.

## 3. Algorithm

The Algorithm for computing maximum [0,2]-factor and characteristic number of graph.

Theorem 2.1 and Theorem 2.2 shows that start from a initiatory [0,2]-factor $\mathbb{R}_1$ of G, if $T_s(\mathbb{R}_1) > 0$ we could to find a P-chain $L_1$ with respect to $\mathbb{R}_1$ then [0,2]-factor $\mathbb{R}_2 = \mathbb{R}_1 \oplus L_1$ have characteristic number $T_S(\mathbb{R}_2) = T_S(\mathbb{R}_1) - 2$ Proceeding in this fashion, until a [0,2]-factor $\mathbb{R}_k$ of G to be found and there can be no P-chain with respect to $\mathbb{R}_k$, then $\mathbb{R}_k$ is the maximum [0,2]-factor of G and $T_S(\mathbb{R}_k)$ is the characteristic number of G. This is done by below Algorithm A.

**Algorithm A:** To Compute a maximum [0,2]-factor of G.

Begin : input a null graph $\mathbb{R}_0$ of G.

$i := 0$

Do

    If Ts($\mathbb{R}_i$)=0 then

    Output a 2-factor $\mathbb{R}_i$ of G ; End

    Else

        Find a P - chain $L_i$ with respect to $\mathbb{R}_i$.

          if there *is* not a P -chain with respect to $\mathbb{R}_i$ then

              Output : A maximal [0,2]-factor $\mathbb{R}_i$ and the characteristic number Ts($\mathbb{R}_i$) of G.:end

          Else

              $\mathbb{R}_{i+1} = \mathbb{R}_i \oplus L_i$

              $i := i+1$

          End if

    End if

Loop

End

In Algorithm A, We must to find a P-chain with respect a [0,2]-factor $\mathbb{R}*$ of G, it can be completed by a DFS. That is, from a vertex $x$ with $d_{\mathbb{R}*}(x) \leq 1$ to find a unmarked adjacent edges $e_1 e_2 ... e_i$ alternately with respected to $E(G) \setminus E(\mathbb{R}*)$ and $E(\mathbb{R}*)$, and mark them, until a end-edge with a end-vertex y and $d_{\mathbb{R}*}(y) \leq 1$ has been found, if for $j < i$ there is no unmarked edge alternating with $e_j$ then back to $e_{j-1}...$ . What is done by below Algorithm.

**Algorithm B:** To find out a P-chain with respect a [0,2]-factor $\mathbb{R}_p$ of G.

Begin: $x_1 \in V(\mathbb{R}_p)$ and $d_{\mathbb{R}_p}(x_1) \leq 1$

$i := 1$

Do

If $i \bmod 2 = 1$ then

  If there is a unmarked edge $e_i = (x_i, x_{i+1}) \in E(G) \setminus E(\mathbb{R}_p)$ then

    $e_i$ to be marked

      If $d_{\mathbb{R}_p}(x_{i+1}) = 1$ and $x_{i+1} \neq x_1$ or $d_{\mathbb{R}_p}(x_{i+1}) = 0$ then

        Output a P-chain: $e_1 e_2 ... e_i$ : end
      Else
        $i := i + 1$
      End if
    Elseif there is no unmarked edge $e_i = (x_i, x_{i+1}) \in E(G) \setminus E(\mathbb{R}_p)$ then
      If $i := 1$ then
        Output: there is no P-chain from vertex $x_1$ : end
      Else
        $i := i - 1$
      End if
    End if
Elseif $i \mod 2 = 0$ then
    If there is a unmarked edge $e_i = (x_i, x_{i+1}) \in E(\mathbb{R}_p)$ then
      $e_i$ to be marked
      $i := i + 1$
    Elseif there is no unmarked edge $e_i = (x_i, x_{i+1}) \in E(\mathbb{R}_p)$ then
      $i := i - 1$
End if
Loop

## The Complexity of Algorithm

    The complexity of all of the algorithm is according to the complexity of the algorithm of finding a P-chain with respect a [0,2]-factor $\mathbb{R}_p$, that is complexity of algorithm B. The algorithm B based on DFS should to search a edge in $E(G)$ once only, so its complexity is $o(|E(G)|)$, therefor, the complexity all of the algorithm is $o(n|E(G)|)$, i.e $o(n^3)$. We can get the below theorem immediately.

    **Theorem 4.1** It can be decision in time $o(n^3)$ that if there is a 2-factor for a graph G, and get a 2-factor when it exist.

## 4   Conclusion

    The chataterristic number of a graph, definited in this paper, is a graph invariant. It shows how far a graph can be [0,2]-factorizations and if and only if which is zero the graph exists 2-factor. The given algorithm make it can be decision and computed out in time $o(n^3)$.

H-cycle is also a 2-factor, therefore a necessary condition of that, a graph is H-graph, is its characteristic number is zero, and it also can be tested in $o(n^3)$.

    The algorithms, given in this paper, has been programmed by VB.

**References**
[1] J. Akiyama and M. Kano. Book of Factors and Factorizations of Graphs. June, 2007. Online version: http://gorogoro.cis.ibaraki.ac.jp/web/papers/Factor    GraphVer1A4.pdf
[2] T. C. Biedl, P. Bose, E. D. Demaine, and A. Lubiw. affecient algorithms for Petersen's matching theorem. *Journal of Algorithms*, 38(1):110{134, 2001.
[3] V. Chvatal. Linear Programming. W. H. Freeman, San Francisco CA, 1983.

[4] P. Diaz-Gutierrez and M. Gopi. Quadrilateral and tetrahedral mesh stripiﬁcation using 2-factor partitioning of the dual graph. *The Visual Computer*, 21(8-10):689{697, 2005.

[5] M. Gopi and D. Eppstein. Single-strip triangulation of manifolds with arbitrary topology. *Computer Graphics*

[6] C. Umans. An algorithm for finding Hamiltonian cycles in grid graphs without holes. *Honors thesis*, Williams College, 1996.

[7] Michael R.Garey and Dvid S.johnson Computers and intractability ,A Guide to the Theory of NP-Completeness. P286,W.H Freeman and Company,San Francisco,1979.